\input amstex\documentstyle{amsppt}  
\pagewidth{12.5cm}\pageheight{19cm}\magnification\magstep1
\topmatter
\title Remarks on affine Springer fibres\endtitle
\author G. Lusztig\endauthor
\address{Department of Mathematics, MIT, Cambridge MA 02139}\endaddress
\thanks{Supported in part by NSF grant DMS-1566618}\endthanks
\endtopmatter   
\document

\define\mpb{\medpagebreak}

\define\si{\sim}

\define\sqc{\sqcup}

\define\lb{\linebreak}

\define\op{\oplus}
   
\define\part{\partial}

\define\m{\mapsto}
\define\do{\dots}

\define\sub{\subset}    

\define\T{\times}
\define\ti{\tilde}
\define\nl{\newline}
\redefine\i{^{-1}}

\define\ov{\overline}
\define\ot{\otimes}

\define\ad{\text{\rm ad}}
\define\Ad{\text{\rm Ad}}

\define\End{\text{\rm End}}

\redefine\b{\beta}

\define\e{\epsilon}

\redefine\o{\omega}
\define\p{\pi}
\define\ph{\phi}
\define\ps{\psi}

\redefine\t{\tau}

\define\x{\xi}

\define\Om{\Omega}

\redefine\L{\Lambda}

\define\CC{\bold C}

\define\EE{\bold E}

\define\NN{\bold N}

\define\ZZ{\bold Z}

\define\cb{\Cal B}

\define\fb{\frak b}

\define\fg{\frak g}

\define\fn{\frak n}

\define\ft{\frak t}

\define\fS{\frak S}

\define\tX{\ti X}

Let $G$ be a simply connected almost simple algebraic
group over $\CC$ and let $\fg$ be the Lie algebra of $G$. Let $B$ be a Borel subgroup of $G$,
let $T$ be a maximal torus of $B$ and let $\ft,\fb$ be the Lie algebras of $T,B$.
Let $\cb$ be the variety of Borel subalgebras of $\fg$. For any nilpotent element
$N\in\fg$ we set $\cb_N=\{\fb\in\cb;N\in\fb\}$ (a Springer fibre). In \cite{KL} an affine 
analogue of $\cb_N$ (``affine Springer
fibre'') was introduced. Let $F=\CC((\e)$, $A=\CC[[\e]]$, where $\e$ is an 
indeterminate and let $\fg(F)=F\ot\fg$ (a Lie algebra over $F$), $L=A\ot\fg$ (a Lie algebra 
over $A$). An element $\x\in\fg(F)$ is said to be topologically nilpotent if
$\lim_{n\to\infty}\ad(\x)^n=0$ in $\End_F(\fg(F))$. Let $\tX$ be the set of all Iwahori 
subalgebras of $\fg(F)$; this is an increasing union of projective varieties over $\CC$. 
According to \cite{KL}, for any regular semisimple, topologically nilpotent element
$\x\in\fg(F)$, the set $\tX_\x=\{I\in\tX;\x\in I\}$
is a nonempty, locally finite union of projective varieties all of the same dimension, say 
$b_\x$. Let $[\tX_\x]$ be the set of irreducible components of $\tX_\x$, a finite or countable
set.

In the remainder of this paper, $h$ denotes a fixed regular element in 
$\ft$. Then $\e h\in\fg(F)$ is regular semisimple, topologically nilpotent so that
the affine Springer fibre $\tX_{\e h}=\{I\in\tX;\e h\in I\}$ is defined. From 
\cite{KL,\S5} we see that $b_{\e h}=\nu$ where $\nu=\dim\cb$.
As in \cite{KL, \S3}, there is a free abelian group $\L$ (see Sec.2) 
of rank equal to the rank of $\fg$ which acts freely on
$\tX_{\e h}$ in such a way that the induced $\L$-action on $[\tX_{\e h}]$ is also free and 
has only finitely many orbits. 
In this paper we will describe a fundamental domain for the $\L$-action on $\tX_{\e h}$.
Namely, let $\fS'$ be the Steinberg variety of triples $(E,\fb_1,\fb_2)$ where 
$\fb_1\in\cb,\fb_2\in\cb$ and $E\in\fb_1\cap\fb_2$ is nilpotent. Let $\fS$ be the fibre at 
$\fb$ of the projection
$\fS'@>>>\cb$, $(E,\fb_1,\fb_2)\m\fb_2$. We can identify $\fS$ with 
$\{(E,\fb_1);\fb_1\in\cb,E\in\fn\cap\fb_1\}$. We state the following result.

\proclaim{Theorem} There is a locally closed subvariety of $\tX_{\e h}$ which is a fundamental
domain for the $\L$-action on $\tX_{\e h}$ such that $\tX_{\e h}$ is isomorphic to $\fS$.
\endproclaim
From the theorem one can deduce some information on the representation of the affine
Weyl group on the vector space $\CC[\tX_{\e h}]$ with basis 
$[\tX_{\e h}]$ defined in
\cite{L2}, see Section 6. 

\mpb

I thank Peng Shan and Zhiwei Yun for discussions.

\subhead 2\endsubhead
Let $U$ be the unipotent radical of $B$. Let $\fn$ be the Lie algebra of $U$. 
Let $G(F),U(F),T(F)$ be the group of $F$-points of $G,U,F$ respectively. Let $G(F)$ be the 
group of $F$-points of $G$. 
Note that $G(F)$ acts naturally
on $\fg(F)$ by the adjoint representation $g:x\m\Ad(g)(x)$.
Let $\L$ be the subgroup of $T(F)$ consisting of the elements $\chi(\e)$ where $\chi$ runs over
the one parameter subgroups $\CC^*@>>>T$ (viewed as homomorphisms $F^*@>>>T(F)$).
Let $X$ be the set of $A$-Lie subalgebras of $\fg(F)$ of the form $\Ad(g)(L)$ for some 
$g\in G(F)$. We shall regard $X$ as an increasing union of projective algebraic varieties
over $\CC$ as in \cite{L1,\S11}. For each $L'\in X$, $L'/\e L'$ inherits from $L'$
a bracket operation and becomes a simple Lie algebra over $\CC$. Let $\p_{L'}:L'@>>>L'/\e L'$
be the obvious map. Let $\cb_{L'}$ be the set of Borel subalgebras of $L'/\e L'$. 
Now $\tX$ consists of all $\CC$-Lie subalgebra of $\fg(F)$ of the form
$\p_{L'}\i(\fb')$ for some $L'\in X$ and some $\fb'\in\cb_{L'}$.
We define $\p:\tX@>>>X$ by $I\m L'$ where $I\sub L'$. Note that $g:I\m\Ad(g)I$
is a well defined action of $G(F)$ on $\tX$ which is transitive. 
According to \cite{KL}, $t:I\m\Ad(t)I$ defines a free action of $\L$ on 
$\tX_{\e h}=\{I\in\tX;\e h\in I\}$ inducing a free
action of $\L$ with finitely many orbits on $[\tX_{\e h}]$.
Let $X_{\e h}=\{L'\in X;\e h\in L'\}$.

If $\x\in\fn(F):=F\ot\fn$ then $\exp(\x)\in U(F)$ is well defined.
Let $\fn(F)'=\op_{i\in\ZZ;i<0}\e^i\fn\sub\fn(F)$.
Let $U(F)'=\{\exp(\x);\x\in\fn(F)'\}\sub U(F)$.
It is well known that any $L'\in X$ can be written in the form
$\Ad(t)\Ad(u)L$ where $t\in\L,u\in U(F)'$ are uniquely determined.
Hence we have a partition $X_{\e h}=\sqc_{t\in\L}X_{\e h,t}$
where $X_{\e h,t}=\{\Ad(t)\Ad(u)L;u\in U(F)',\e h\in\Ad(u)L\}$ is a locally closed subset of 
$X_{\e h}$. Let $\tX_{\e h,t}=\p\i(X_{\e h,t})$. This is a 
locally closed subset of $\tX_{\e h}$.
 Let $\Om=X_{\e h,1},\ti\Om=\tX_{\e h,1}=\p\i(\Om)$. Note that 

(a) $\tX_{\e h}=\sqc_{t\in\L}\Ad(t)\ti\Om$
\nl
as a set. Thus, $\ti\Om$ is a fundamental domain for the $\L$-action on $\tX_{\e h}$.
Let $\o=\{\EE\in\fn(F)';\Ad(\exp(\EE))(\e h)\in L\}$. In preparation for the proof of the
theorem we will prove the following result.

\proclaim{Lemma 3}The map $\EE=\e\i E_1+\e^{-2}E_2+\e^{-3}E_3+\do@>>>E_1$ is a bijection
$\ph:\o@>\si>>\fn$. (Here $E_1,E_2,E_3,\do$ is a sequence of elements of $\fn$ with $E_i=0$ for
large $i$.)
\endproclaim
The equation defining $\o$ is $\exp(\ad(\EE))(\e h)\in L$ that is
$$\align&\e h+\sum_{i\ge1}\e^{-i+1}[E_i,h]+(1/2)\sum_{i,j\ge1}\e^{-i-j+1}[E_i,[E_j,h]]\\&+
(1/6)\sum_{i,j,k\ge1}\e^{-i-j-k+1}[E_i,[E_j,[E_k,h]]]+\do\in L,\endalign$$
that is 
$$\align&\sum_{i\ge2}\e^{-i+1}[E_i,h]+(1/2)\sum_{i,j\ge1}\e^{-i-j+1}[E_i,[E_j,h]]\\&+
(1/6)\sum_{i,j,k\ge1}\e^{-i-j-k+1}[E_i,[E_j,[E_k,h]]]+\do\in L,\endalign$$
that is
$$\align&[E_r,h]=-(1/2)\sum_{i,j\ge1,i+j=r}[E_i,[E_j,h]]\\&-
(1/6)\sum_{i,j,k\ge1,i+j+k=r}[E_i,[E_j,[E_k,h]]]+\do\tag a\endalign$$
for $r=2,3,\do$.
In the right hand side we have $i<r,j<r,k<r$, etc. Hence if $E_{r'}$ is known for 
$r'<r$ then $[E_r,h]$ is a well defined element of $\fn$. Hence $E_r$ is a well defined 
element of $\fn$. (Note that $E\m[E,h]$ is a vector space isomorphism $\fn@>\si>>\fn$.)

It remains to show that $E_r=0$ for large $r$. 
For $r\ge1$ let $\fn^r$ be the subspace of $\fn$ spanned by all iterated brackets of $r$
elements of $\fn$. (Thus, $\fn^1=\fn$, $\fn^2$ is spanned by $[a,b]$ with $a,b$ in $\fn$,
$\fn^3$ is spanned by $[[a,b],c]]$ with $a,b,c$ in $\fn$, etc.) Note that 

(b) $E\m[E,h]$ is an isomorphism $\fn^r@>>>\fn^r$ for any $r\ge1$.
\nl
We show by induction on $r$ that
$$E_r\in \fn^r\text{ for }r=1,2,\do\tag c$$
For $r=1$ this is clear. Assume now that $r\ge2$. From (a) and the induction hypothesis
we deduce that $[E_r,h]\in\fn^r$. Using (b) we see that for some $E'\in\fn^r$ we have
$[E_r,h]=[E',h]$, hence $[E_r-E',h]=0$, hence $E_r=E'$. Thus $E_r\in\fn^r$, proving (c).
Since $\fn^r=0$ for large $r$ we see that $E_r=0$ for large $r$.
This completes the proof of the lemma.

\subhead 4\endsubhead
For $E\in\fn$ we set $u_E=\exp(\EE)\in U(F)'$ where $\EE=\ph\i(E)$ (see Lemma 3). Note that 
$\Ad(u_E)(\e h)\in L$. Now $\EE\m\Ad(\exp(-\EE))L$ is a bijection $\ps:\o@>\si>>\O$. Hence 
$\ps':=\ps\ph\i:\fn@>>>\Om$ is a bijection. We have $\ps'(E)=\Ad(u_E\i)L$. We show:

(a) {\it Let $E\in\fn$ and let $L_E=\Ad(u_E\i)L\in X$. Note that $\e h\in L_E$. Then
$\p_{L_E}(\e h)\in L_E/\e L_E$ and $\p_L(-[E,h])\in L/\e L$ correspond to each other under the
Lie algebra isomorphism $\t_E:L/\e L@>\si>>L_E/\e L_E$ induced by $\Ad(u_E\i):L@>\si>>L_E$.}
\nl
We must show that $Ad(u_E)(\e h)=-[E,h]\mod\e L$ or that $\Ad(\exp(\EE))(\e h)=-[E,h]\mod\e L$
where $\EE=\e\i E_1+\e^{-2}E_2+\e^{-3}E_3+\do$ corresponds to $E=E_1$ as in Lemma 3.
Thus we must show that
$$\align&\e h+\sum_{i\ge1}\e^{-i+1}[E_i,h]+(1/2)\sum_{i,j\ge1}\e^{-i-j+1}[E_i,[E_j,h]]\\&+
(1/6)\sum_{i,j,k\ge1}\e^{-i-j-k+1}[E_i,[E_j,[E_k,h]]]+\do =-[E_1,h]\mod\e L,\endalign$$
or that
$$\align&\sum_{i\ge2}\e^{-i+1}[E_i,h]+(1/2)\sum_{i,j\ge1}\e^{-i-j+1}[E_i,[E_j,h]]\\&+
(1/6)\sum_{i,j,k\ge1}\e^{-i-j-k+1}[E_i,[E_j,[E_k,h]]]+\do\in\e L.\endalign$$
But the left hand side is actually zero, by the proof of Lemma 3. This proves (a).

\mpb

From (a) we deduce:

(b) {\it the map $\b\m\t_E(\b)$ is a bijection 
$\{\b\in\cb_L;\p_L(-[E,h])\in\b\}@>>>\{\b'\in\cb_{L_E};\p_{L_E}(\e h)\in\b'\}$.}
\nl
Taking union over all $E\in\fn$ and using the bijection
$\ps':\fn@>>>\Om$ we deduce 

(c) {\it the map $(E,\b)\m\p\i_{L_E}(\t_E(\b))$ is a bijection 
$\{(E,\b)\in\fn\T\cb_L;\p_L(-[E,h])\in\b\}@>\si>>\ti\Om$.}
\nl
We consider the bijection

(d) $\{(E,\b)\in\fn\T\cb_L;\p_L(-[E,h])\in\b\}@>>>\fS$ 
\nl
given by $(E,\b)\m(-[E,h],\fb_1)$ where $\fb_1\in\cb$ is defined by $\p_L(\fb_1)=\b$.
The composition of the inverse of (d) with the bijection (c) is a bijection 

(e) $\fS@>\si>>\ti\Om$.
\nl
From the proof we see that the bijection (e) is an isomorphism of algebraic varieties.
This proves the theorem.

\subhead 5\endsubhead
Let $NT$ be the normalizer of $T$ in $G$ and let $W=NT/T$ be the Weyl group. For any $w\in W$ 
let $\cb_w$ be the variety consisting of all $\fb_1\in\cb$ such that $(\fb,\fb_1)$
are in relative position $w$. Note that $\cb_w$ is isomorphic to $\CC^{|w|}$ where $|w|\in\NN$ 
is the length of $w$. Let $\fS_w=\{(E,\fb_1)\in\fS;\fb_1\in\cb_w\}$. The second projection
$\fS_w@>>>\cb_w$ makes $\fS_w$ into a vector bundle with fibres of dimension $\nu-|w|$.
Hence $\fS_w$ is isomorphic to $\CC^\nu$ as an algebraic variety. We 
have a partition $\fS=\sqc_{w\in W}\fS_w$ (as a set) with $\fS_w$ locally closed in $\fS$
(the closure of $\fS_w$ in $\fS$ is denoted by $\ov{\fS_w}$).
Hence we have a partition $\ti\Om=\sqc_{w\in W}\ti\Om_w$ (as a set) where $\ti\Om_w$ 
corresponds to $\fS_w$ under 4(e). Note that $\ti\Om_w$ is isomorphic to $\CC^\nu$ as an 
algebraic variety and that $\ti\Om_w$ is locally closed in $\ti\Om$.
For $w\in W,t\in\L$ we set $\ti\Om_{w,t}=\Ad(t)\ti\Om_w$. Using 2(a) we see that

(a) $\tX_{\e h}=\sqc_{(w,t)\in W\T\L}\ti\Om_{w,t}$
\nl
as a set, where $\ti\Om_{w,t}$ is locally closed in $\tX_{\e h}$ and is isomorphic to $\CC^\nu$.
Let $\ov{\ti\Om_{w,t}}$ be the closure of $\ti\Om_{w,t}$ in $\tX_{\e h}$. Note that
$\ti\Om_{w,t}$ is open dense in $\ov{\ti\Om_{w,t}}$. Since $\tX_{\e h}$ is of pure dimension 
$\nu$, we see that 

(b) {\it $(w,t)\m\ov{\ti\Om_{w,t}}$ is a bijection $W\T\L@>\si>>[\tX_{\e h}]$.}
\nl
In particular,

(c) {\it the number of $\L$-orbits on $[\tX_{\e h}]$ is equal to the order of $W$.}
\nl
A result closely related to (c) (but not (c) itself) appears in \cite{TS}.

\subhead 6\endsubhead
Let $[\fS]$ be the set of irreducible components of $\fS$ (a finite set naturally indexed by
$W$ by $w\m\ov{\fS_w}$). The bijection 5(b) gives rise to an imbedding
$[\fS]@>>>[\tX_{\e h}]$, $\ov{\fS_w}\m\ov{\ti\Om_{w,1}}$ hence to an imbedding of vector spaces

(a) $\CC[\fS]@>>>\CC[\tX_{\e h}]$
\nl
with bases $[\fS],[\tX_{\e h}]$. Springer has shown that $W$ acts naturally on 
$\CC[\fS]$ (this is known to be the 
regular representation of $W$ in a nonstandard basis). In \cite{L2} it is shown that the affine Weyl group of $G$ acts naturally on
$\CC[\tX_{\e h}]$. Hence, by restriction, $W$ acts on $\CC[\tX_{\e h}]$. From the definitions
we see that the imbedding (a) is compatible with the $W$-actions.

\enddocument